\documentclass[11pt]{article}
\usepackage{amsmath,amsthm,amsfonts,amssymb,graphicx,graphpap}
\usepackage{hyperref}
\usepackage[all]{xy}

\title{Fixed Subgroups of Endomorphisms of Free Products}
\author{Mihalis Sykiotis }
\makeatletter
\newtheorem{thm}{Theorem}[section]
\newtheorem{lem}[thm]{Lemma}
\newtheorem{cor}[thm]{Corollary}

\theoremstyle{definition}

\theoremstyle{remark}

\def\classification{\@ifnextchar [{\@xfootnotetext}%
{\begingroup\let\protect\noexpand\xdef\@thefnmark{}
\endgroup\@footnotetext}}
\begin{document}

\classification {2000 {\sl Mathematics Subject Classification.}
20E06, 20E36. } \maketitle
\begin{abstract} Let $G=\ast_{i=1}^{n}G_{i}$ and let $\phi$ be a symmetric
endomorphism of $G$. If $\phi$ is a monomorphism or if $G$ is a
finitely generated residually finite group, then the fixed
subgroup $Fix(\phi)=\{g\in G:\phi(g)=g\}$ of $\phi$ has Kurosh
rank at most $n$.
\end{abstract}

\section{Introduction}
In \cite{BH}, Bestvina and Handel proved the Scott conjecture,
which says that if $\phi$ is an automorphism of a free group of
rank $n$, then the subgroup $Fix(\phi)$ of elements fixed by
$\phi$ has rank at most $n$. Their result was generalized by
several authors in various directions. See, for example, \cite{IT,
CT, DV, Ber, Sy}. In particular, the result of Bestvina and Handel
was generalized both to arbitrary endomorphisms of free groups by
Imrich and Turner \cite{IT} and to automorphisms of free products
by Collins and Turner \cite{CT}.

In this note, following the main idea of \cite{IT}, we show that
in many interesting cases the study of fixed subgroups of
endomorphisms of free products is reduced to that of
automorphisms, thereby obtaining new generalizations of
Bestvina-Hadel's result.

\section{Preliminaries}
Let $G=\ast_{i=1}^{n}G_{i}$ and let $H$ be a non-trivial subgroup
of $G$. By the Kurosh subgroup theorem, $H$ is a free product
$H=\ast_{i\in I}H_{i}\ast F$, where $F$ is a free group and every
factor $H_{i}$ is the intersection of $H$ with a conjugate of a
free factor $G_{i}$. In the case where the rank $r(F)$ of $F$ and
the cardinality $|I|$ of $I$ (which may be empty) are finite, the
{\sl Kurosh rank of H with respect to the given splitting of $G$}
is defined to be the sum $r(F)+|I|$. We will usually omit the
phrase ``with respect ... splitting of $G$", when the splitting of
$G$ is clear from the context.
\par
Following \cite{IT}, given a group $G$ and an endomorphism $\phi$
of $G$, we define the {\sl stable image} $\phi^{\infty}(G)$ of
$\phi$ to be the intersection $\cap_{n=1}^{\infty}\phi^{n}(G)$.
Clearly $\phi^{\infty}(G)$ is invariant under $\phi$ and contains
$Fix(\phi)$. Thus $Fix(\phi)=Fix(\phi_{\infty})$, where
$\phi_{\infty}:\phi^{\infty}(G)\rightarrow \phi^{\infty}(G)$
denotes the restriction of $\phi$ to $\phi^{\infty}(G)$. The key
observation is that if $\phi$ is a monomorphism, then
$\phi_{\infty}$ is an automorphism. To see this, let $g\in
\phi^{\infty}(G)$ be any element. Then for every $n$ there exists
an element $g_{n}$ of $\phi^{n}(G)$ such that $g=\phi(g_{n})$.
Since $\phi$ is injective, $g_{1}=g_{n}$ for all $n$ and hence
$g_{1}\in \phi^{\infty}(G)$. This gives surjectivity of
$\phi_{\infty}$.
\par
Now, the basic idea can be described briefly as follows. Suppose
that $G$ is a free product and that $\phi_{\infty}$ is an
automorphism sending non-infinite-cyclic factors of the stable
image onto conjugates of themselves. By \cite[Theorem 6.12]{Sy},
the Kurosh-rank of $Fix(\phi_{\infty})$ does not exceed the Kurosh
rank of $\phi^{\infty}(G)$. Thus to find an upper bound for the
Kurosh rank of $Fix(\phi)$, we need to know something about the
kurosh rank of $\phi^{\infty}(G)$. By \cite[Theorem 6.5]{Sy1}, the
Kurosh rank of $\phi^{\infty}(G)$ is bounded above by the maximum
of the Kurosh ranks of the images $\phi^{n}(G)$. In the case where
$G$ is a free group of rank $n$, it is immediate that the rank of
every image $\phi^{n}(G)$ is less than or equal to $n$ while in
the case of a free product is not. However, we will see that this
happens in many cases, in which we obtain that the Kurosh rank of
$Fix(\phi)$ does not exceed the Kurosh rank of $G$.

\section{Main Results}

We start with the following result which has been obtained
independently by Swarup ~\cite{Sw}.
\begin{lem}\label{lem:1} Let $G=\ast_{i=1}^{n}G_{i}\ast F$ and $H=\ast_{j=1}^{m}H_{j}\ast
F'$, where each factor $G_{i}$ is not infinite cyclic and $F$,
$F'$ are free groups, and let $\phi:G\rightarrow H$ be an
epimorphism such that each factor $G_{i}$ is mapped by $\phi$ into
a conjugate of some $H_{j}$. Then $n+r(F)\geq m+r(F')$ and
$r(F)\geq r(F')$.
\end{lem}

\begin{proof} By renumbering if necessary, we can assume that
$H_{1},\dots,H_{m_{0}}$, $m_{0}\leq m$ are the factors of $H$
whose conjugates contain the non-trivial images of $G_{i}$,
$i=1,\dots,n$ under $\phi$. Note that $m_{0}\leq n$. If $N$ and
$K$ are the normal subgroups of $G$ and $H$ generated by
$G_{i},i=1,\dots,n$ and $H_{j},j=1,\dots,m_{0}$ respectively, then
$\phi(N)\subseteq K$, and so $\phi$ induces an epimorphism
$\Phi:F\cong G/N \rightarrow H/K \cong H_{m_{0}+1}\ast \cdots \ast
H_{m}\ast F'$. It follows that $r(F)\geq d( H_{m_{0}+1}\ast \cdots
\ast H_{m}\ast F')= d(H_{m_{0}+1})+ \cdots +d(H_{m})+r(F')\geq
m-m_{0}+r(F')$, and the lemma follows.
\end{proof}

Let $G=\ast_{i=1}^{n}G_{i}$ and $H=\ast_{i=1}^{m}H_{i}$. A
homomorphism $\phi:G\rightarrow H$ is said to be
\textsl{symmetric} if each non-infinite-cyclic free factor of $G$
is mapped by $\phi$ into a conjugate of some non-infinite-cyclic
free factor of $H$. For example, if each factor $G_{i}$ is freely
indecomposable, then each injective homomorphism is symmetric.

The next lemma shows that symmetric automorphisms of free products
map non-infinite-cyclic factors onto conjugates of themselves and
therefore ~\cite[Theorem 6.12]{Sy} can be applied.
\begin{lem}\label{lem:2} Let $G=\ast_{i=1}^{n}G_{i}\ast F$ and let $\phi$ be an automorphism of $G$. If each factor
$G_{i}$ is mapped by $\phi$ into a conjugate of some $G_{j}$, then
$G_{i}$ is mapped by $\phi$ onto this conjugate.
\end{lem}

\begin{proof}
Suppose on the contrary that there is a factor, say $G_{1}$, such
that $\phi(G_{1})$ is properly contained in $gG_{i_{1}}g^{-1}$,
$g\in G$. By ~\cite[Theorem 7]{St} there is a free product
decomposition $G=\ast_{i=1}^{n}G'_{i}\ast F'$ of $G$ such that
$\phi(G'_{i})=G_{i}$, $i=1,\dots,n$ and $\phi(F')=F$. If
$x=\phi^{-1}(g)$, then $\phi(x^{-1}G_{1}x)\subset
G_{i_{1}}=\phi(G'_{i_{1}})$ and thus $x^{-1}G_{1}x\subset
G'_{i_{1}}$. Since $G'_{i_{1}}$ properly contains $x^{-1}G_{1}x$,
there is a free product decomposition
$G'_{i_{1}}=x_{1}G_{1}x_{1}^{-1}\ast K$, obtained from the initial
decomposition of $G$, where $K$ is a non-trivial subgroup of $G$.
Thus $G=G_{1}\ast\cdots \ast G_{n}\ast F=\ast_{i\neq i_{1}}
G'_{i}\ast x_{1}G_{1}x_{1}^{-1}\ast K\ast F'$.
\par
Now we consider the map $\psi:\{1,\dots,n\}\rightarrow
\{1,\dots,n\}$, defined as follows: $\psi(i)=j$ if and only if
$\phi(G_{i})$ is contained in a conjugate of $G_{j}$. The
injectivity of $\phi$ implies that $\psi$ is well-defined, while
the proof of Lemma \ref{lem:1} shows that $\psi$ is surjective,
and hence bijective. We conclude that the normal subgroup $N$ of
$G$ generated by $G_{1},\dots,G_{n}$ is contained in the normal
subgroup $N'$ of $G$ generated by $G_{1},G'_{i}$, $i\neq i_{1}$.
Thus we have an epimorphism $\Phi:G/N\cong F\rightarrow G/N'\cong
K\ast F'$. Since $K$ is non-trivial, it follows that $r(F)>r(F')$,
which contradicts the fact that the groups $F$ and $F'$ are
isomorphic.
\end{proof}

\begin{thm}\label{thm:1} Let $G=\ast_{i=1}^{n}G_{i}\ast F$, where each factor
$G_{i}$ is not infinite cyclic and $F$ is a free group. If $\phi:
G\rightarrow G$ is a symmetric monomorphism of $G$, then the fixed
subgroup $Fix(\phi)$ of $\phi$ has Kurosh rank at most $n+r(F)$.
\end{thm}

\begin{proof}
By the remarks preceding Lemma \ref{lem:1}, it suffices to show
that the stable image $\phi^{\infty}(G)$ of $\phi$ has Kurosh rank
at most $n+r(F)$, and that the automorphism $\phi_{\infty}$ of
$\phi^{\infty}(G)$ is symmetric, since the theorem is true for
symmetric automorphisms of free products ~\cite{Sy}.
\par
First, we note that for each $k\geq 0$, the epimorphism
$\phi_{k}:\phi^{k}(G)\rightarrow \phi^{k+1}(G)$ obtained by
restricting $\phi$ to $\phi^{k}(G)$ is symmetric (where
$\phi^{0}(G)=G$), which implies that $\phi^{k}(G)$ has Kurosh rank
at most $n+r(F)$ for all $k$ by Lemma ~\ref{lem:1}. To see this,
let $\phi^{k}(G)\cap xG_{i}x^{-1}$ be a non-infinite-cyclic free
factor of $\phi^{k}(G)$ (with respect to the free product
decomposition of $\phi^{k}(G)$ inherited from this one of $G$).
The assumption that $\phi$ is symmetric implies that there is an
index $j(i)\in \{1,\dots,n\}$ and an element $g_{i}\in G$ such
that $\phi(G_{i})\subseteq g_{i}G_{j(i)}g_{i}^{-1}$. Thus
$\phi(\phi^{k}(G)\cap xG_{i}x^{-1})\subseteq \phi^{k+1}(G)\cap
\phi(x)g_{i}G_{j(i)}g_{i}^{-1}\phi(x)^{-1}$. The latter group is a
subgroup of $\phi^{k+1}(G)$ which stabilizes a vertex in any
$G$-tree constructed from the given free product decomposition of
$G$. It follows that $\phi(\phi^{k}(G)\cap xG_{i}x^{-1})$ is
contained in a $\phi^{k+1}(G)$-conjugate of a free factor of
$\phi^{k+1}(G)$ and hence $\phi_{k}$ is symmetric. The same
argument shows that the automorphism $\phi_{\infty}$ is symmetric
as well.
\par
Since each term of the decreasing sequence of subgroups
\[G\supseteq \phi(G)\supseteq \phi^{2}(G)\supseteq \cdots
\supseteq  \phi^{k}(G)\supseteq \cdots\] has Kurosh rank at most
$n+r(F)$, ~\cite[Theorem 6.5]{Sy1} implies that the stable image
$\phi^{\infty}(G)$ of $\phi$ also has Kurosh rank at most
$n+r(F)$.
\end{proof}

\begin{cor} Let $\phi$ be a monomorphism of a free product
$\ast_{i=1}^{n}G_{i}$ of freely indecomposable groups. Then
$Fix(\phi)$ has Kurosh rank at most $n$.
\end{cor}

In view of the preceding theorem, it is natural to seek conditions
under which a free product endomorphism becomes ``finally" a
monomorphism. The second proof of the Hopficity of finitely
generated residually finite groups sketched in ~\cite{Ha},
actually shows that the restriction of an endomorphism of a
residually finite group to its stable image is a monomorphism (see
also ~\cite[Lemma 1]{Hi}). For completeness, we include the
argument here.
\begin{lem}\label{lem:3} Let $\phi$ be an endomorphism of a finitely generated residually finite group
$G$. Then the restriction
$\phi_{\infty}:\phi^{\infty}(G)\rightarrow \phi^{\infty}(G)$ of
$\phi$ to $\phi^{\infty}(G)$ is a monomorphism.
\end{lem}
\begin{proof}
Let $1\neq g\in ker(\phi_{\infty})$. Then $\phi(g)=1$ and for each
positive integer $n$ there is $g_{n}\in G$ such that
$g=\phi^{n}(g_{n})$. Since $G$ is residually finite there is a
finite group $\Gamma$ and a homomorphism $\pi:G\rightarrow \Gamma$
with $\pi(g)\neq 1$. We consider the sequence of homomorphisms
$\pi_{n}=\pi\circ\phi^{n}:G\rightarrow \Gamma$. Then $1\neq
\pi(g)=\pi\big(\phi^{n}(g_{n})\big)=\pi_{n}(g_{n})$. On the other
hand,
$\pi_{m}(g_{n})=\pi\big(\phi^{m}(g_{n})\big)=\pi\big(\phi^{m-n}(g)\big)=1$
whenever $m>n$. It follows that there are infinitely many distinct
homomorphisms from the finitely generated group $G$ to the finite
group $\Gamma$, a contradiction.
\end{proof}

\begin{thm}\label{thm:2} Let $G=\ast_{i=1}^{n}G_{i}\ast F$ be a finitely generated
residually finite group, where each factor $G_{i}$ is not infinite
cyclic and $F$ is a free group. If $\phi$ is a symmetric
endomorphism of $G$, then the fixed subgroup $Fix(\phi)$ of $\phi$
has Kurosh rank at most $n+r(F)$.
\end{thm}
\begin{proof}  The arguments of Theorem \ref{thm:1} show that the stable image $\phi^{\infty}(G)$ of
$\phi$ has Kurosh rank at most $n+r(F)$ and that the restriction
$\phi_{\infty}$ of $\phi$ to $\phi^{\infty}(G)$ is a symmetric
endomorphism. By Lemma \ref{lem:3}, $\phi_{\infty}$ is a
monomorphism, so Theorem \ref{thm:1} applies.
\end{proof}

\begin{thm} Let $G=\ast_{i=1}^{n}G_{i}$ be a free product of
finitely generated nilpotent and finite groups. If $\phi$ is an
endomorphism of $G$, then the fixed subgroup $Fix(\phi)$ of $\phi$
has Kurosh rank at most $n$.
\end{thm}
\begin{proof}
Since each quotient of a nilpotent group is freely indecomposable,
each of the epimorphisms $\phi_{k}:\phi^{k}(G)\rightarrow
\phi^{k+1}(G)$ satisfies the hypothesis of Lemma \ref{lem:1}. This
implies that $\phi^{\infty}(G)$ has Kurosh rank at most $n$. By
Lemma \ref{lem:3}, $\phi_{\infty}$ is a monomorphism. Also, it is
easy to see that $\phi_{\infty}$ is symmetric. The theorem now
follows by Theorem \ref{thm:1}.
\end{proof}

Department of Mathematics and Statistics, University of Cyprus,
P.O. Box 20537, 1678 Nicosia, Cyprus

{\sl E-mail address:} msikiot@ucy.ac.cy, msykiot@math.uoa.gr

\end{document}